\documentclass[12pt]{article}       
\usepackage{a4}
\usepackage{amsmath,amssymb,amsthm} 
\usepackage[all]{xy}
\newcommand{\N}{\mbox{$\mathbb N$}}        
\newcommand{\Q}{\mbox{$\mathbb Q$}}        

\newcommand{\sq}{\mbox{$\displaystyle\sqcup\!\sqcup$}}

\theoremstyle{plain}               
\newtheorem{defn}{Definition}[section]
\newtheorem{lemma}[defn]{Lemma}
\newtheorem{proposition}[defn]{Proposition}
\newtheorem{corollary}[defn]{Corollary}

\theoremstyle{definition}          
\newtheorem{remark}[defn]{Remark}
\newtheorem{example}[defn]{Example}

\begin{document}

\title{Planar Shuffle Product, Co-Addition and the
non-associative Exponential}
\author{by L.~Gerritzen}
\date{}

\maketitle

\begin{abstract}

In this note we introduce the concept of a shuffle product $\sq$
for planar tree polynomials and give a formula to compute the
planar shuffle product $S \ \sq \ T$ of two finite planar reduced
rooted trees $S, T.$

It is shown that $\sq$ is dual to the co-addition $\Delta$ which
leads to a formula for the coefficients of $\Delta(f).$

It is also proved that $\Delta(EXP) = EXP {\hat \otimes} EXP$
where $EXP$ is the generic planar tree exponential series, see
[G]. Systems of quadratic relations for the coefficients of EXP
are derived.
\end{abstract}


{\bf \large{Introduction}}

\medskip

Let $K \{x \}_\infty$ be the algebra of planar tree polynomials
over a field $K$where the $m$-ary product is obtained by the
$K$-multi-linear extension of the $m$-ary planar grafting of
finite planar reduced rooted trees. There is a unique $K$-algebra
homomorphism

$$\Delta : K \{x\}_\infty \to K \{x\}_\infty \otimes K
\{x\}_\infty$$

mapping $x$ onto $x \otimes 1 + 1 \otimes x$ where $x$ denotesd
the tree with a single vertex.

 This coproduct is also called co-addition.
Dually there is a shuffle product

$$\sq : K \{x \}_\infty \otimes K \{x\}_\infty \to K
\{x\}_\infty$$

which is commutative and associative.

It can be defined by the contradiction of planar trees onto sets
of leaves.

For any $f \in K \{x\}_\infty$ the coefficients $(v,w(\Delta(f))$
of $\Delta(f)$ can be computed by a formula involving the
coefficients of the shuffle product  $V \ \sq \ W,$ see Prop. 4.3.

The co-addition extends continuously to a $K$-algebra homomorphism

$$\hat{\Delta} : K \{\{x\}\}_\infty \to K \{\{x\}\}_\infty
\hat{\otimes} K \{\{x\}\}_\infty$$

of formal power series completions of $K\{x\}_\infty$ and
$K\{x\}_\infty \otimes K\{x\}_\infty.$

It is shown in Section 5 that

$$\hat {\Delta}(EXP) = EXP \hat{\otimes} EXP$$

if $EXP$ is the generic non-associative exponential series, see
[G], Section 3. One obtain systems of quadratic relations for the
coeffcients of $EXP$ in which the coefficients of planar shuffle
products appear.

This result can be used to extend the canonical projection of [R],
Chap. 3 or [L], to the non-associative setting. This application
has been suggested to me by J. L. Loday.

\bigskip\bigskip

\begin{section}{Reductions of planar rooted trees}

Let $T$ be a finite planar rooted tree, see [G]. For any vertex
$a$ of $T$ we denote by $val_T(a)$ the valence of $a$ in $T$ which
is the number of edges of $T$ incident with $a.$

Let $$ar_T(a) := \left \{\begin{array}{cc} val_T (a)\ \ \ \ \ :&
\hbox{ if }a \hbox{ is the root of } T\\ val_T(a) - 1  :& \hbox{
otherwise }
\end{array}\right.$$

It is called the arity of the vertex $a$ in $T.$

\begin{defn} $T$ is called reduced, if $ar_T(a) \not=1$ for all
vertices $a$ of $T.$
\end{defn}

\medskip

In the following construction we associate a reduced planar rooted
tree Red$(T)$ to a given finite planar rooted tree $T$ which might
not be reduced.

For vertices $a, b$ of $T$ we denote by $[a,b]_T$ the smallest
connected subgraph of $T$ which contains $a$ and $b.$ It is a line
and consists of the vertices and edges of the simple path between
$a$ and $b.$

\medskip

Let $$R^0 := \{ a \in T^0= \hbox{ set of vertices of } T:ar_T(a)
\not= 1 \}$$
 and $$\bar{R} := \{\{a,b\} : a,b \in R^0, a \not=b, [a,b]^0_T
\cap R^0 = \{a, b\}\}.$$

Let $\rho_R$ be the vertex of $R^0$ closest to the root $\rho_T$
of $T.$

\begin{proposition}
$(R^0, \bar{R}, \rho_R)$ is a planar reduced rooted tree; it is
called the reduction of $T$ and is also denoted by Red$(T).$
\end{proposition}

This result will follow from Proposition 1.4 below.

\begin{remark} If the arity of the root $\rho_T$ of a planar
rooted tree $T$ is 1, then Red$(T)$ = Red$(T-\rho_T)$ where
$T-\rho_T$ is obtained from $T$ by deleting $\rho_T$ and the
unique edge in $T$ incident with $\rho_T.$
\end{remark}

\begin{proposition}

\begin{itemize}
\item[(i)] If $T_1,...,T_m$ are finite planar reduced rooted trees and
$$T = \cdot_m(T_1,...,T_m)$$ is the grafting of $T_1,...,T_m,$
then $T$ is a planar reduced rooted tree if $m \ge 2.$
\item[(ii)] If $T_1,...,T_m$ are finite planar rooted trees and $m \ge 2,$ then

$$Red(\cdot_m(T_1,...,T_m)) = \cdot_m( Red(T_1),..., Red(T_m)).$$

\end{itemize}
\end{proposition}

\medskip
\begin{proof}
\begin{itemize}
\item[1)] By definition the set $T^0$ of vertices of $T$ consists of
the root $\rho_T$ and the disjoint union the set $T^0_i$ of
vertices of $T_i$ for $1 \le i \le m.$

Obviously for any vertex of $T^0_i$ we get

$$ar_{T_i}(a) = ar_T(a)$$

while $ar_T(\rho_T) = m.$

This shows statement (i).

\item[2)] Let $T = \cdot_m(T_1,...,T_m).$ Then $R^0 = \{a \in T^0
ar_T(a) \not= 1 \}$ is the disjoint union $R^0_1,...,R^0_m$ where
$R^0_i = \{a \in T^0_i: ar_{T_i} (a) \not= 1 \}$ as the arity of
the root of $T$ is equal to $m \ge 2.$ Any line $[a,b]_T$ in $T$
for which all vertices $\not= a, \not=b$ have arity 1 in $T$ is
therefore lying in a subtree $T_i$ of $T.$ From this the statement
(ii) follows.
\end{itemize}
\end{proof}

\begin{example} The reduction $T$ of the following planar
rooted trees $T_1, T_2, T_3$ are all equal to $T = 
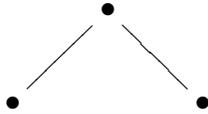
\xymatrix
{&&\bullet\ar@{-}[dl]\ar@{-}[dr]&\\ &\bullet&&\bullet}$

\medskip

In all cases the root of the trees is the most upward vertex.

\bigskip\bigskip


$
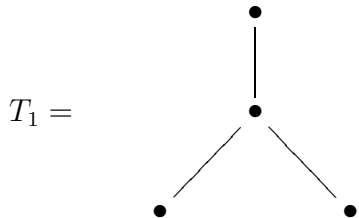
\xymatrix{&&\bullet\ar@{-}[d]&&\\ T_1 =
&&\bullet\ar@{-}[dl]\ar@{-}[dr]&&\\ &\bullet &&\bullet&&}$

\bigskip\bigskip

 $
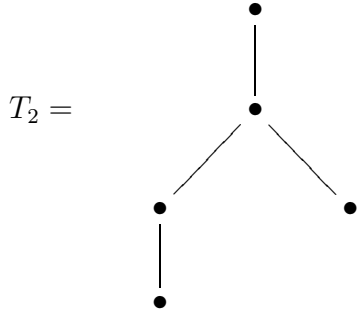
\xymatrix{ &&\bullet\ar@{-}[d]&&\\
 T_2=&&\bullet\ar@{-}[dl]\ar@{-}[dr]&&\\
 &\bullet\ar@{-}[d] &&\bullet\\
  &\bullet&&}$

  \bigskip\bigskip

$
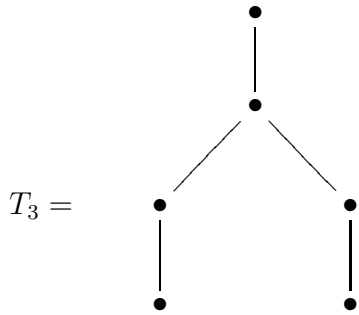
\xymatrix{ &&\bullet\ar@{-}[d]&&\\
 &&\bullet\ar@{-}[dl]\ar@{-}[dr]&&\\
T_3= &\bullet\ar@{-}[d] &&\bullet
  \ar@{-}[d]\\ &\bullet&& \bullet&}$

\bigskip\bigskip

 \end{example}
\end{section}

\begin{section}{Contractions}

Denote by ${\bf PRT}$ the set of isomorphism classes of finite
planar reduced rooted trees.

Let $T \in {\bf PRT}$ and $I$ be a subset of the set $L(T) := \{a
\in T^0: ar_T(a) = 0\}$ of leaves of $T.$

Then $T_I := \bigcup_{a \in I} [a, \rho_T]_T$ is a connected
subgraph of $T.$ It is a subtree of $T$ because every connected
subgraph of a tree is again a tree. The tree $T_I$ is the smallest
subtree of $T$ containing $I$ and $\rho_T.$ We consider $T_I$ as a
rooted tree with root $\rho_T.$ Then $L(T_I) = I$ because all
vertices $b$ of $[a, \rho_T]_T$ different from $a$ have arity 1 in
$[a, \rho_T]_T.$

\begin{defn}

The contraction $T \vert I$ of $I$ onto $I$ is defined to be the
reduction $Red(T_I)$ of $T_I.$
\end{defn}

\begin{proposition}

Let $T \in{\bf PRT}, ar(T) = m \ge 2 \ T = T_1 \cdot T_2 \cdot ...
\cdot T_m$ with $T_i \in {\bf PRT}.$

Let $I \subset L(T) = L(T_1) \dot\cup ... \dot\cup L(T_m)$ and
$I_j= I \cap L(T_j).$

Then $T\vert I = T_1 \vert I_1 \cdot ... \cdot T \vert I_m.$
\end{proposition}

\begin{proof}

By standard considerations.
\end{proof}
\end{section}

\begin{section}{Shuffle product}

Let $K$ be a field and ${\bf A} = K \{x\}_\infty$ the $K$-algebra
with unit of planar tree polynomials over $K.$

A $K$-base of ${\bf A}$ is given by ${\bf PRT'} = \{1\}\dot\cup
{\bf PRT}$ where 1 denotes the empty treeand the $m$-ary
multiplication $\cdot_m : {\bf A}^m \to{\bf A}$ is the
$K$-multi-linear extension of the $m$-ary grafting on ${\bf
PRT'},$ see [G].

For $S, T, V \in {\bf PRT'}$ let $N_{S,T}(V) := \{I \subset L(V) :
V \vert I = S, V \vert I^c = T \}$ where $I^c := L(V) - I$ is the
complement of $I$ in $L(V).$

If $deg(V) \not= deg(S) + deg(T),$ then $N_{S,T}(V)$ is empty.

\begin{defn}

$$S \ \sq \ T := \sum_{V \in {\bf PRT'}} \sharp N_{S,T}(V) \cdot
V$$

is called the planar shuffle product of $S$ and $T.$
\end{defn}

\begin{remark}

This planar shuffle product is a generalization of the usual
shuffle product as defined in [R], section (1.4).
\end{remark}
\begin{proposition}

There is a unique $K$-bilinear map

$$\sq : K \{x\}_\infty \times K \{x\}_\infty \to K \{x\}_\infty$$

such that $\sq \ (S,T)$ is the shuffle product of $S$ and $T$ for
any $S,T \in {\bf PRT'}.$

It is called the shuffle product on $K \{x\}_\infty.$

Moreover $\sq$ is commutative and associative.
\end{proposition}

\begin{proof}
\begin{itemize}
\item[1)] The existence and uniqueness of the $K$-bilinear map
follows immediately from the fact that ${\bf PRT'}$ is a $K$-base
of $K\{x\}_\infty.$

\item[2)] The commutativity of $\sq$ follows because $N_{S,T}(V) =
N_{T,S}(V).$ A bijection between $N_{S,T}(V)$ and $N_{T,S}(V)$ is
given by $I \to I^c.$

\item[3)] Let now $R, S, T \in{\bf PRT'}$ and $N_{R,S,T}(V) :=
\{(I_1, I_2, I_3) : I_j \subseteq L(V) \ I_1, I_2, I_3$ is a
disjoint union of $L(V), V \vert _1, R, \ V \vert I_2 = S, \
V\vert I_3 = T \}.$ It is an easy exercise to show that $$R \ \sq
\ S) \ \sq \ T = \sum_{W \in {\bf PRT'}} \sharp N_{R,S,T}(V) \cdot
V.$$
\end{itemize}
The right hand side is also equal to $R \ \sq (S \ \sq \ T).$
\end{proof}
\end{section}

For $k \in \N_{\ge 1}$ denote by {\underline k} the set $\{1, 2,
... , k \}$ of integers between 1 and $k.$

For any non-empty finite subset $\alpha$ ogf $\N_{\ge 1}$ denote
by $\alpha [i]$ the i-th element of $\alpha$ relative to the
natural order on $\alpha$ for $i \in{\underline k}, \ k := \sharp
\alpha.$ Thus the map

$$i \mapsto \alpha [i]$$

is a bijection between ${\underline k}$ and $\alpha$ and $\alpha$
[i] $< \alpha[i + 1]$ for $1 \le i < k.$

We represent $\alpha$ by the word $\alpha [1]\cdot  \alpha [2]
\cdot ... \cdot \alpha[k].$

For any $k \in \N_{\ge 2}$ denote by $\Gamma_k :=\{(\alpha,
\beta): \alpha \subseteq {\underline k}, \beta \subseteq
{\underline k}, \alpha \not= \emptyset \not= \beta, \alpha \cup
\beta = {\underline k}\}$ the set of pairs $(\alpha, \beta)$ of
non-empty subsets of ${\underline k}$ whose union is the whole set
${\underline k}.$

For any $m, n \in \N_{\ge 1}$ denote $\Gamma_k(m, n) := \{(\alpha,
\beta) \in \Gamma_k : \sharp \alpha = m, \sharp \beta = n \}.$
Obviously the map $(\alpha, \beta) \mapsto (\beta, \alpha)$ is a
bijection between $\Gamma_k(m,n)$ and $\Gamma_k(n, m).$

\begin{lemma}
\begin{itemize}
\item[i)] $\Gamma_k(m, n) \not= \empty$ if and only if $m \le
k, n \le k, k \le m+n.$
\item[ii)] If $m \le k, n \le k, k \le m + n,$ then $\sharp
\Gamma_k(m, n) = \frac{k!}{(k-n)!(k-m)!(m+n-k)!}$
\end{itemize}
 \end{lemma}

 \begin{proof}
 \begin{itemize}
 \item[1)] (i) is immediate, as for $(\alpha, \beta) \in
 \Gamma_k(m,n)$ one has $\alpha \subseteq {\underline k}, \beta
 \subseteq {\underline k}, \alpha \cup \beta = {\underline k}.$
 from which $\sharp \alpha \le k, \sharp \beta \le k, \sharp
 \alpha + \sharp \beta \ge \sharp {\underline k} = k$ follows.
 \item[2)] Assume the assumption of (ii). There are ${k \choose
 n}$ subsets $\alpha$ of ${\underline k}$ with $n$ elements. There
 are ${n \choose m+n-k}$ subsets $\gamma$ of $\alpha$ with $m + n -
 k$ elements. Let $\beta := ({\underline k} - \alpha) \cup \gamma.$
Any pair of $\Gamma_k(m,n)$ is obtained by this procedure. Thus
$\sharp \Gamma_k(m, n) = {k \choose n} \cdot {n \choose m+n-k}.$
\end{itemize}
 \end{proof}

\begin{example}

\medskip

$ \Gamma_2(1,1) = \{(1, 2), (2,1)\}$

$ \Gamma_2(1,2) = \{(1,1 2), (2,12)\}$

$ \Gamma_2(2,1) = \{(12, 1), (12,2)\}$

$ \Gamma_2(2,2) = \{(12,1 2) \}$

 $ \Gamma_2(m,n) = \emptyset$ if $m$ or $n$ is $\ge 3$

\medskip

$\Gamma_k(1,n)  \not= \emptyset$ if and only if $n \le k \le n +
1.$

$\Gamma_n(1,n) = \{(i, {\underline n}) : i \in {\underline n}$

$\Gamma_{n+1}(1,n) = \{(i, {\underline {n+1}} - \{i\}) : i \in
{\underline n}\}.$

\medskip

Especially

\medskip

$\Gamma_3(1,3) = \{(1, {\underline 3}), (2, {\underline 3}), (3,
{\underline 3})\}$

$\Gamma_4(1,4) = \{(1, 234), (2, 134), (3, 124), (4,123)\}.$
\end{example}

\bigskip

Let ${\bar \Gamma}_k(m,n) := \Gamma_k(1,1) \cup \Gamma_k(m,1) \cup
\Gamma_k(1, n) \cup\Gamma_k(m,n).$

It is the set of all pairs $(\alpha, \beta)$ of subsets of
${\underline k}$ such that $(\sharp \alpha$ is 1 or $m, \  \sharp
\beta$ is 1 or $n$ and $\alpha \cup \beta = {\underline k}.$

Let $S, T \subset {\bf PRT}, m = ar(S), n = ar(T), S = S_1 \cdot
S_2 \cdot ... \cdot S_m, T = T_1 \cdot T_2 \cdot .. \cdot T_n$
with $S_i, T_j \in {\bf PRT}.$

For any $(\alpha, \beta) \in {\bar \Gamma}_k(m, n)$ we are going
to define a polynomial $R_{(\alpha, \beta)} \in K \{x\}_\infty.$
One has to consider several cases:

\medskip

 Case 1:

\begin{itemize}
\item[] $\sharp \alpha = \sharp \beta = 1.$ Then $(\alpha, \beta)
\in \Gamma_2(1,1)$ and $R_{(\alpha, \beta)}:=  S\cdot T \hbox{ if
} \alpha = 1, \beta = 2  and R_{\alpha, \beta} := T \cdot S \hbox{
if } \alpha = 2, \beta = 1.$

\end{itemize}

\medskip

Case 2:

\begin{itemize}
\item[] $\sharp \alpha = 1, \sharp \beta = m > 1, \alpha =
\{\nu\}.$

Then $(\alpha, \beta) \in \Gamma_k(1, n)$ and $k$ is either $n$ or
$n+1.$

If $k = n + 1,$ then $\beta = {\underline{n+1}} - \{\nu\}$ and
$$R_{\alpha, \beta)} := T_1 \cdot T_2 \cdot ... \cdot
T_{\nu-1}\cdot S \cdot T_\nu \cdot T_{\nu+1} \cdot T_n$$ which is
a tree of arity $n+1.$

If $k = n,$ then $\beta = {\underline n}$ and $R_{(\alpha, \beta)}
:= T_1 \cdot T_2 \cdot ... \cdot T_{\nu-1}\cdot (T_\nu \ \sq \ S)
\cdot T_{\nu+1}\cdot ... \cdot T_n$ which is a polynomial of arity
$n.$
\end{itemize}

\medskip

Case $\prime$:

\begin{itemize}
\item[] $\sharp \beta = 1, \beta = \{\mu\}, \sharp \alpha = n >
1.$ Then $(\alpha, \beta) \in \Gamma_k(m, 1)$ and $k = m$ or $m +
1.$

If $k = m + 1,$ then $\alpha = {\underline {m+1}} - \{\mu\}$ and

$R_{(\alpha, \beta)} := S_1 \cdot S_2 \cdot ... \cdot S_{\mu-1}
\cdot T \cdot S_\mu \cdot S_{\mu+1}\cdot ... \cdot S_m$ which is a
tree of arity $m + 1.$ If $k = m,$ then $\alpha = {\underline m}$
and $R_{(\alpha, \beta)};= S_1 \cdot S_2 \cdot ... \cdot
S_{\mu-1}\cdot (S_\mu \ \sq \ T) \cdot S_{\mu+1} \cdot ... \cdot
S_m$ which is a polynomial of arity $m.$
\end{itemize}

\medskip

Case 3:
\begin{itemize}
\item[] $\sharp \alpha = m > 1, \sharp \beta = n > 1.$ Then
$R_{(\alpha, \beta)}:= P_1 \cdot P_2 \cdot ... \cdot P_k$

where $P_j = S_i$ if $j \in \alpha, j \notin \beta$ and $\alpha
[i] = j$

where $P_j = T_i$ if $j \in \beta, j \notin \alpha$ and $\beta[i]
= j$

and where $P_j = S_i \ \sq \ T_l$ if $j \in \alpha \cap \beta$ and
$\alpha[i] = j, \beta [l] = j.$
\end{itemize}

\medskip

\begin{proposition}
$$S \ \sq \ T = \sum^{n+m}_{k = 2} \sum_{(\alpha, \beta) \in {\bar
\Gamma}_k(m,n)} R_{(\alpha, \beta)}.$$
\end{proposition}

\medskip

\begin{proof}

\begin{itemize}
\item[1)] Let $N_{S,T} = \displaystyle \cup_{V \in {\bf PRT}} N_
{S,T} (V).$

Its elements can be seen as pairs $(V, I)$ with $V \in {\bf PRT},$
$I \subset L(V)$ and $I \in N_{S, T}(V).$ Obviously $S \ \sq \ T =
\sum V$ where the summation is extended over $N_{S,T}.$
\item[2)] Let $(V, I) \in N_{S, T}$ and $V$ be of arity $k.$ Then
$V = V_1 \cdot V_2 \cdot ...  \cdot V_k$ with $V_i \in {\bf PRT}$
and $L(V) = L(V_1) \dot\cup \cdot ... \cdot \dot \cup L(V_k).$

The type $ (\alpha, \beta)$ of $(V, I)$ is defined as follows:

$\alpha := \{ i \in {\underline k} : I \cap L_i(V) \not=
\emptyset\}$

$\beta:= \{i \in {\underline k} : I^c \cap L_i(V) \not=
\emptyset\}$

where $I^c$ is the complement of $I$ in $L(V).$

It is easy to see that $(\alpha, \beta) \in {\bar \Gamma}_k(m,
n).$ If $\alpha \cap \beta = \emptyset$ it is obvious that $V =
R_{(\alpha, \beta)}$ and $I = \displaystyle \cup_{i \in
\alpha}L_i(V).$

\item[3)] If $\alpha \cap \beta \not= \emptyset$ then the situation
is more complicated to describe.

What happens can already be seen if $m = n = k$ and $\alpha =
\beta = {\underline k}.$

Then $R_{(\alpha, \beta)} = (S_1 \ \sq \ T_1) \dot ... \cdot /S_n
\ \sq \ T_n) =  \sum V_1 \cdot V_2 \cdot ... \cdot V_n$ where the
summation is extended over the set of n-tupels $((V_1, I_1), ...,
(V_n,I_n))$ with $(V_i, I_i) \in N_{S_i, T_i}$ where $S = S_1
\cdot S_2 ... S_n, T = T_1 \cdot T_2 \cdot ... \cdot T_n.$

This is obviously equal to the sum $\displaystyle \sum_{(V,I) \in
N_{S,T}} V,$ where $N'_{S,T}$ is the subset of $N_{S,T}$
consisting of all $(V, I)$ for which the type is equal to
$({\underline n, \underline n}.$

\item[4)] The general case follows by a combination of the arguments
 in 2) and 3).
\end{itemize}
\end{proof}
\begin{example} $x \ \sq \ x = 2x^2.$

If $S = x,$ and $n \ge 2,$ then

$x \ \sq \ T = x \cdot T + t \cdot x + \displaystyle
\sum^n_{i=1}T_1 \cdot ... \cdot T_{i-1}\cdot x \cdot T_i \cdot
T_{i+1}\cdot ... \cdot T_n + \displaystyle \sum^n_{i=1}T_1 \cdot
... \cdot T_{i-1} \cdot (x \ \sq \ T_i) \cdot T_{i+1}\cdot ...
\cdot T_n.$

It follows from the Proposition because ${\bar \Gamma}_k(1,n) =
\Gamma_k(1,1) \cup \Gamma_k(1,n)$ and from the explicit
description of $\Gamma_k(1,n)$ above.
\end{example}

\begin{section}{Co-Addition}

Let ${\bf A} = K\{x\}_\infty$ and $B = {\bf A} \otimes_K {\bf A}.$

\begin{proposition}

For any $m \in \N_{\ge 2}$ there is a unique $K$-multi-linear map

$$\cdot_m : {\bf B}^m \to {\bf B}$$

such that $\cdot_m(f_1 \otimes g_1,...,f_m \otimes g_m) = (f_1
\cdot f_2 \cdot ... \cdot f_m) \otimes (g_1 \cdot g_2 \cdot ...
\cdot g_m)$ for all $f_i, g_i \in {\bf A}.$

${\bf B}$ together with this string $(\cdot_m)_{m \ge 2}$ of
operations is considered as an algebra.
\end{proposition}
\begin{proof}
As in the classical case of one binary multiplication.
\end{proof}

\begin{proposition}

There is a unique unital $K$-algebra homomorphism

$$\Delta : {\bf A} \to {\bf A} \otimes{\bf  A}$$

such that $\Delta(x) = x \otimes 1 + 1 \otimes x.$

$\Delta$ is called co-addition or coproduct induced by addition.
\end{proposition}

\begin{proof}
By standard methods.
\end{proof}

For $T \in {\bf PRT'}$ let

$$c_T : {\bf A} \to K $$

be the $K$-linear map with $$c_T(S) = \delta_{ST} =\left\{
\begin{array}{cc} 1 & : S=T\\
0 & :S\not=T\end{array}\right. \\$$ for any $S \in {\bf PRT'}.$

\medskip

It is called the coefficient map relative to $T.$

For $V,W \in {\bf PRT'}$ let

$$c_{V,W} : {\bf A} \otimes {\bf A} \to K $$

be the $K$-linear map such that
 $$c_{V,W}(S \otimes T) = \left\{ \begin{array}{cc} 1 :& (V,W) = (S,T)\\
 0 :& \hbox{ otherwise } \end{array}\right.\\$$

for $S,T \in {\bf PRT'}.$ It is called the coefficient map
relative to $V \otimes W.$

\medskip

The following formula for planar polynomials generalizes the
classical formula, see [R].
\begin{proposition}

For any $f \in K \{x \}_\infty$ and $V,W \in {\bf PRT'}$ one gets

$$c_{V,W}(\Delta(f)) = \sum_{T \in {\bf PRT'}} c_T(f) \cdot c_T(V
\ \sq \ W).$$
\end{proposition}
\begin{proof}
\begin{itemize}
\item[1)] For any $T \in {\bf PRT'}$ we obtain $\Delta(T) =
\displaystyle \sum_{I cL(T)} (T \vert I) \otimes (T \vert I^c)$ as
can be seen by induction on $deg(T),$ see also [GH]. It follows
that $c_{V,W}(\Delta(T)) = c_T(V  \ \sq \ W).$ Then the formula is
proved if $f$ is a monomial.
\item[2)] Both sides of the formula are $K$-linear and thus they
are equal because of step 1).
\end{itemize}
\end{proof}

%
%
%
%
%

\end{section}

\begin{section}{Co-Addition of $EXP$}

Let ${\bf \hat A}$ be the $x$-adic completion of ${\bf A} = K
\{x\}_\infty$ and ${\bf \hat B}$ be the $(x \otimes 1, 1 \otimes
x)$-adic completion of ${\bf B}.$ Then $\Delta :{\bf  A} \to {\bf
A} \otimes {\bf A}$ extends to a continuous coproduct
$$\hat{\Delta} : {\bf \hat A} \to {\bf \hat B} = {\bf \hat
A}\hat{\otimes} {\bf \hat A}.$$

\begin{proposition}

Let $EXP$ be the generic exponential in $K \{\{x\}\}_\infty,$ see
{\rm [G]}, $K = \Q(q)$ field of rational functions in $q$ over
$\Q.$

Then

$$\hat{\Delta}(EXP) = EXP \hat{\otimes} EXP$$
\end{proposition}

\begin{proof}
\begin{itemize}
\item[1)] Let $k \in \N_{\ge 2}$ and $exp_k(x)$ be the $k$-ary
exponential series in $\Q \{\{x\}\}_\infty.$ Let $f =
\hat{\Delta}(exp_k(x)).$

We want to show that $f = exp_k(x)\hat \otimes exp_k(x).$
\item[2)] For $k \in {\bf \hat B}$ let $ord_{\bf B}(k) := min\{deg(V) +
deg(W) : c_{V,W}(h) \not= 0, V,W \in {\bf PRT'}\}.$ Obviously
$ord_{{bf \hat B}} \Delta(g) = ord_{{\bf \hat A}}(g)$ for $g \in
{\bf \hat A}$ and thus

$$ord_{{\bf \hat B}}(f - (1 \otimes 1 + x \otimes 1 + 1 \otimes
x)) \ge 2$$
\item[3)] There is a continuous $K$-automorphism
$$\varphi_k : {{\bf \hat B}}\to{\bf  \hat B}$$
 such that

$$\varphi_k(x \hat{\otimes} 1) = k(x \hat{\otimes} 1)$$
$$\varphi_k(1 \hat {\otimes} x) = k(1 \hat{\otimes} x).$$

Fix $k$ and let $f = \Delta(EXP^k) = \Delta(EXP(kx)) =
\varphi_k(f).$

Denote by $f_r$ the homogeneous part of $f$ degree $r.$

Then
 $$ \sum_{V,W  \in {\bf PRT'}\atop
degV+degW=r}c_{V,W}(f) \cdot V \otimes W.$$
 As $\varphi_k(f_r) =k^rf_r$we get from the functional equation
 for $f$ that

$$(f^k) = k^rf_r$$

if $(f^k)_r$ is the homogeneous part of $f^k$ of degree $r$ which
is

$$\sum_{degV+degW=r} c_{V,W}(f^k) \cdot V \otimes W.$$

As in [G] the Prop. (3.1) one can show that $f$ must be equal to
$exp_k \hat{\otimes}exp_k.$
\item[4)]Let $R$ be the subring of all rational function $f$ in
$\Q(q)$ whose poles are roots of unity $\not= 1.$

Let $\pi_k : R \to \Q$ be the $\Q$-algebra homomorphism with
$\pi_k(q) = k$ for $k \in \N_{\ge 2}.$

Now $\pi_k$ induces algebra homomorphisms $$\hat{\pi}_k: R
\{\{x\}\}_\infty \to \Q\{\{x\}\},$$ $$ \hat{\pi}^2_k : R
\{\{x\}\}_\infty \hat{\otimes} R \{\{x \}\}_\infty  \to \Q
\{\{x\}\}_\infty \hat{\otimes} \Q\{\{x\}\}_\infty$$ and

$$\hat{\pi}^2_k(\Delta(k)) = \Delta(\hat{\pi}_k(k))$$ for $k  \in
R \{\{x\}\}_\infty.$

As $\hat{\pi}_k(EXP) = exp_k$ and $\hat{\pi}^{(2)}(\Delta(EXP)) =
exp_k \hat{\otimes}exp_k$ we get $$\Delta(EXP) = EXP \hat{\otimes}
EXP.$$
\end{itemize}
\end{proof}

Let $a(T) = c_T(EXP)$ be the coefficient of $EXP$ relative to $T
\in {\bf PRT'}.$ Recall that $a(1) = a(x) = 1$ and the recursive
relation

$$a(T) = \frac{{q \choose m}}{q^n-q} a(T_1) \dot ... \cdot
a(T_m)$$

if $m = ar(T), T = T_1 \cdot ... \cdot T_m$ and $n = deg(T),$ see
[G].

In the following Corollary one obtains a quadratic relation for
the system $(a_T)_{T \in {\bf PRT'}}.$

\begin{corollary}
$a(V) \cdot a(W) = \displaystyle \sum_T a_T c_T(V \ \sq \ W).$

\end{corollary}

\begin{example}

Let $m \ge 3$ and $S = x \ \sq \ x^m.$ One can show that

$S = (m + 1)x^{m+1} + 2(x^2 \cdot x \cdot ... \cdot x + x \cdot
x^2 \cdot x \cdot ... \cdot x + ... + x \cdot x \cdot ... \cdot x
\cdot x^2) + x \cdot x^m + x^m \cdot x$ by applying Example (3.6).

From the recursive formula for the coefficients of $EXP$  one gets

$$a(x^m) = \frac{{q\choose m}}{q^m-q}$$

$$a(x^2 \cdot x \cdot ... \cdot x) = \frac{{q\choose
m}}{q^{m+1}-q} \cdot \frac{1}{2} = a(x \cdot ... \cdot x \cdot x^2
\cdot x \cdot ... \cdot x)$$ if the factor $x$ occurs $m$ times,
as $a(x^2) = \frac{1}{2}.$

$$a(x^m \cdot x) = a(x \cdot x^m) = \frac{{q\choose 2}}{q^{m+1}-q}
\cdot \frac{{q\choose m}}{q^m-q}.$$

The right hand side in the formula of Corollary (5.2) is $$b =
\displaystyle \sum_{T\in {\bf PRT'}} c_T(S) \cdot a(T).$$

We are substituting the values above for $a(T)$ and obtain

$$b = (m+1)\frac{{q\choose m+1}}{q^{m+1}-q} + 2m \frac{{q \choose
m}}{q^{m+1}-q} \cdot1/2 + 2 \frac{{q\choose 2}}{q^{m+1}-q}\frac
{{q\choose m}}{q^m-q} $$ as $deg((x^{m+1}) = m + 1, \ deg(x^2
\cdot x \cdot ... \cdot x) = m + 1$ if the factor $x$ occurs $m$
times and $deg(x \cdot x^m) = deg(x^m \cdot x) = m + 1.$

Thus $$b =\frac{{q \choose m}}{q^m - q} ((m+1)
\frac{(\frac{q-m}{m+1}\cdot (q^m-q))}{(q^{m+1} -q)}+ m \cdot
\frac{q^m-q}{q^{m+1}-q}+ \frac{q(q-1)}{q^{m+1}-q}) = \frac{{q
\choose m}}{q^m-q}$$ because

$$(q-m) (q^m - q) + m(q^m-q) + q(q-1) = q^{m+1} - q$$ and $b =
a(x) \cdot a(x^m)= a(x^m) = \frac{{q \choose m}}{q^m-q}$ which is
in accordance with Corollary (5.2).

\end{example}
\end{section}


\end{document}